\documentclass[12pt,reqno]{amsart}

\usepackage[T1]{fontenc}
\usepackage{ae,aecompl}
\usepackage{graphics,cancel,graphicx}
\usepackage{epsfig,psfrag,manfnt}
\usepackage{mathrsfs,amsmath}
\usepackage{graphicx,mathabx}
\usepackage{bbm,subfigure,rotating,bm,float,mathdots,wasysym,scalerel,color,xcolor}
\usepackage{mathrsfs,setspace}
\usepackage{graphics,cancel,graphicx,mathabx}
\usepackage{epsfig,psfrag,manfnt}
\usepackage{bbm,subfigure,rotating,bm,float,mathdots,wasysym,scalerel}
\usepackage{array}
\usepackage{comment}

\definecolor{cerulean}{rgb}{0.0, 0.48, 0.65}
\definecolor{burntorange}{rgb}{0.8, 0.33, 0.0}
\definecolor{green}{rgb}{0.3,0.6,0.3}
\definecolor{red}{rgb}{1, 0.2, 0.1}
\usepackage[T1]{fontenc}
\usepackage{ae,aecompl}
\usepackage{mathrsfs,setspace}
\usepackage{graphics,cancel,graphicx,mathabx}
\usepackage{epsfig,psfrag,manfnt}
\usepackage{bbm,subfigure,rotating,bm,float,mathdots,wasysym,scalerel,yfonts}
\usepackage{array}

\usepackage{relsize}
\usepackage{exscale,xparse}

\usepackage{tikz-cd}
\usepackage{tikz}
\usetikzlibrary{positioning} 

\usepackage{array}
\usepackage[hidelinks]{hyperref}
\usepackage{xcolor}
\hypersetup{
    colorlinks,
    linkcolor={burntorange!50!black},
    citecolor={blue!50!black},
    urlcolor={blue!80!black}
}

\usepackage{relsize}
\usepackage{exscale,xparse}
\usepackage{tikz-cd}
\usepackage{tikz}
\usetikzlibrary{positioning} 
\usepackage{array}
\usepackage{hyperref}

\usepackage[all,cmtip]{xy}
\usepackage{amssymb,amsmath,amsfonts,amsthm,color}

\usepackage{enumitem}
\setlist[itemize]{align=parleft,left=0pt..1em}
\setlist[enumerate]{align=parleft,left=0pt..1.63em}

\newtheorem{theorem}{Theorem}[section]
\newtheorem{lemma}[theorem]{Lemma}
\newtheorem{corollary}[theorem]{Corollary}
\newtheorem{proposition}[theorem]{Proposition}

\theoremstyle{definition}

\theoremstyle{definition}

\theoremstyle{definition}

\theoremstyle{definition}
\newtheorem{example}[theorem]{Example}

\begin{document}

\setlength{\abovedisplayskip}{3pt}
\setlength{\belowdisplayskip}{3pt}



\title[]{On the faithful flatness \\
of some modules arising in analysis}

\author[]{Amol Sasane}
\address{Mathematics Department, London School of Economics, 
Houghton Street, London WC2A 2AE, UK}
\email{A.J.Sasane@LSE.ac.uk}
 
\maketitle

\vspace{-0.51cm}

\begin{abstract}The notion of faithful flatness of a module over a commutative ring is studied 
 for two $R$-modules $M$ arising in functional analysis, 
where $R$ is a Banach algebra and $M$ is a Hilbert space. The following results are shown:  

\smallskip

\noindent 
 If $X$ is a locally compact Hausdorff topological space, and $\mu$ is a positive Radon measure on $X$, 
 then $L^2(X,\mu)$ is a flat $L^\infty(X,\mu)$-module. Moreover:
 
 \noindent ${\scriptstyle\bullet}$  If $\mu$ is $\sigma$-finite, 
 then for every finitely generated, nonzero, proper \phantom{${\scriptstyle\bullet}$ }ideal $\mathfrak{n}$ of $L^\infty(X,\mu)$, there holds $\mathfrak{n}L^2(X,\mu)\subsetneq L^2(X,\mu)$. 
 
 \noindent ${\scriptstyle\bullet}$ 
 If $X$ is the union of an increasing family of Borel sets $U_n$, $n\!\in\! {\mathbb{N}}$, 
 \phantom{${\scriptstyle\bullet}$ }such that for each $n\in {\mathbb{N}}$, $\overline{U_n}$ is compact and $\mu(U_{n+1}\setminus U_n)>0$, \phantom{${\scriptstyle\bullet}$ }then 
$L^2(X,\mu)$ is not a faithfully flat $L^\infty(X,\mu)$-module. 

\smallskip

\noindent 
In addition, it is shown that the classical Hardy space $H^2$ is a flat, but not a faithfully flat $H^\infty$-module, which answers a 2005 question of Alban Quadrat.
\end{abstract}

\section{Introduction and main results.}

\noindent The aim of this article is to investigate the algebraic notion of faithful flatness 
for two natural modules arising in analysis.

\subsection{The notions of flatness and faithful flatness}

The notion of flatness was introduced by J-P. Serre in \cite{Ser}. 
Recall that a left $R$-module $M$ over a ring $R$ is {\em flat} 
if for every injective linear module morphism $\varphi : K\to L$  
of right $R$-modules $K$ and $L$, the map
$$
\varphi\otimes_R M: K\otimes_R M\to L\otimes_R M
$$
is also injective, where $\varphi \otimes _{R}M$ is the map 
induced by $k \otimes m \mapsto \varphi ( k ) \otimes  m$. 
If $R$ is a unital commutative ring, then an $R$-module $M$ 
is  flat if and only if  for every linear relation 
$$
\textstyle 
\sum\limits_{i=1}^n
r_im_i=0
$$
(where $n\in {\mathbb{N}}$, 
$m_1,\cdots, m_n\in M$ and $r_1,\cdots, r_n\in R$), 
there exist
\begin{itemize}
\item $k\in {\mathbb{N}}$,  
\item $\rho_{ij}\in R$, $i\in \{1,\cdots, n\}$, $j\in \{1,\cdots, k\}$, and 
\item $\mu_1,\cdots, \mu_n\in M$ 
\end{itemize}
such that we have 
\begin{itemize}
\item $\sum\limits_{i=1}^n
r_i \rho_{ij}=0$ for all $ j\in \{1,\cdots, k\}$, and  
\item  $
 m_i=\textstyle 
\sum\limits_{j=1}^k 
\rho_{ij} \mu_j
$ for each $i\in \{1,\cdots,n\}$.
\end{itemize}
An $R$-module $F$ is {\em faithfully flat}  if it is flat and for every
nonzero $R$-module $M$, $F \otimes_R M\neq 0$. Equivalently, 
an $R$-module $F$  is faithfully flat if it is flat and for every maximal ideal $\mathfrak{m}$ of $R$, $\mathfrak{m} F  \neq  F$, 
where 
$$
\textstyle 
\mathfrak{m}F:=\{\sum\limits_{i=1}^n r_i f_i: n\in \mathbb{N}, \text{ and for all }i\in \{1,\cdots,n\}, \;r_i\in \mathfrak{m} \text{ and } f_i\in F\}.
$$
 Of course $\mathfrak{m}F\subset  F$. 
(Throughout, we use the notation $A\subset B$ to mean inclusion of the set $A$ in $B$, allowing $A$ and $B$ to be equal, and we use $A\subsetneq B$ for  strict inclusion.) 

\subsection{Motivation} 

Besides the interest in the notion of faithful flatness of modules from the commutative algebra perspective, 
the faithful flatness of particular $R$-modules $M$, for some concrete Hilbert spaces $M$ and Banach algebras $R$ arising in control theory, play a key role in the stabilisation problem for linear control systems, see \cite{Qua}. 
In fact, one of the results established in this article, namely that $H^2$ is not a faithfully flat $H^\infty$-module (see page~\pageref{pageref_sarason}), answers a question from 2005 raised in \cite[Remark~1]{Qua}. 

\subsection{First main result} 

Let $X$ be a locally compact Hausdorff topological space, and $\mu$ be a positive Radon measure on $X$. 
Let $L^\infty(X,\mu)$ be the set of all $f:X\to {\mathbb{C}}$ such that 
$$
\textstyle
\begin{array}{rcl}
\|f\|_\infty &=&{\text{ess$\;\!$sup}}_{x\in X} |f(x)|\\[0.1cm]
 &=&\inf\{ u\ge 0: |f(x)|\le u \text{ for }\mu\text{-almost all }x\in X\}.
 \end{array}
$$
We identify functions in $L^\infty(X,\mu)$ that differ on any set with $\mu$-measure equal to $0$. 
With pointwise operations and the norm $\|\cdot\|_\infty$, $L^\infty(X,\mu)$ is a Banach algebra. 
Let $L^2(X,\mu)$ be the set of all $f:X\to {\mathbb{C}}$ such that 
$$
\textstyle 
\|f\|_2^2=\int_X |f(x)|^2 d\mu(x)<\infty.
$$
We identify functions in $L^2(X,\mu)$ that differ on any set with $\mu$-measure equal to $0$. 
Then $L^2(X,\mu)$ is a Hilbert space with pointwise operations and the norm $\|\cdot\|_2$ (which is induced by an inner product). 
It is clear that with the action of $L^\infty(X,\mu)$ on $L^2(X,\mu)$ given by pointwise multiplication, $L^2(X,\mu)$ is an $L^\infty(X,\mu)$-module. We show that $L^2(X,\mu)$ is a flat $L^\infty(X,\mu)$-module (Proposition~\ref{18_9_24_1810}). 

Moreover, we have the following results: 
 \begin{itemize}
\item If $\mu$ is $\sigma$-finite, 
 then for every finitely generated, nonzero, proper ideal $\mathfrak{n}$ of $L^\infty(X,\mu)$, there holds $\mathfrak{n}L^2(X,\mu)\subsetneq L^2(X,\mu)$ (Proposition~\ref{19_9_24_1616}). 
 \item
 If $X$ is the union of an increasing family of Borel sets $U_n$, $n\!\in\! {\mathbb{N}}$, such that for each $n\in {\mathbb{N}}$, $\overline{U_n}$ is compact and $\mu(U_{n+1}\setminus U_n)>0$, then 
$L^2(X,\mu)$ is not a faithfully flat $L^\infty(X,\mu)$-module. This is an immediate consequence of   Theorem~\ref{19_9_24_1617}. 
\end{itemize}

\vspace{-0.51cm}

\subsection{Second main result} 
 Let 
 $
 {\mathbb{D}}\!=\!\{z\in {\mathbb{C}}\!:\!|z|\!<\!1\}$ 
 be the unit disc. Also, let  
 ${\mathbb{T}}\!=\!\{z\in {\mathbb{C}}\!:\! |z|\!=\!1\}$,  and $
 \overline{{\mathbb{D}}}\!=\!{\mathbb{D}}\cup {\mathbb{T}}$.
   The set of complex-valued holomorphic functions on ${\mathbb{D}}$ is denoted by ${\mathcal{O}}({\mathbb{D}})$. 
 Let the {\em Hardy algebra} $H^\infty$ be the Banach algebra consisting of all bounded  functions $f\!\in\! {\mathcal{O}}({\mathbb{D}})$, equipped with pointwise operations, and the supremum norm, given by  $\|f\|_\infty=\sup_{z\in {\mathbb{D}}}|f(z)|$ for $f\in H^\infty$. If $f\in H^\infty$, then the radial limits 
$$
\textstyle 
f(e^{i\theta}):=\lim\limits_{r\to 1} f(re^{i\theta})
$$
exist for almost all $\theta \in (-\pi,\pi]$, and define a  function in $L^\infty({\mathbb{T}})$, which we have also denoted by $f$ above. 
The {\em Hardy Hilbert space} $H^2$ is the set of all $h\in {\mathcal{O}}({\mathbb{D}})$ such that 
$$
\textstyle
\|h\|_2^2:=\sum\limits_{n=0}^\infty  |h_n|^2<\infty,
\;\text{ where }\;
h(z)=\sum\limits_{n=0}^\infty  h_n z^n \;\; (z\in {\mathbb{D}}).
$$
Then $H^2$ is a Hilbert space with the inner product corresponding to the norm $\|\cdot\|_2$ defined above. 
It can be shown that 
$$
\textstyle 
\|h\|_2=\sup\limits_{r\in (0,1)} \big(\frac{1}{2\pi} \int_{-\pi}^{\pi} |h(re^{i\theta})|^2 \;\!d\theta\big)^{\frac{1}{2}}.
$$
If $f\in H^2$, then again the radial limits 
$$
\textstyle 
h(e^{i\theta}):=\lim\limits_{r\to 1} h(re^{i\theta})
$$
exist for almost all $\theta \in (-\pi,\pi]$, and now define a boundary function in $L^2({\mathbb{T}})$, which we have also denoted by $h$.  The set of these boundary functions coincides with the set of all $h\in L^2({\mathbb{T}})$ whose negatively indexed Fourier coefficients vanish, that is, $\hat{h}_{-n}=0$ for all $n\in {\mathbb{N}}$, where
$$
\textstyle 
\hat{h}_{m}:=\frac{1}{2\pi}\int_{-\pi}^{\pi} h(e^{i\theta})e^{-im\theta} \;\! d\theta \quad (m\in {\mathbb{Z}}).
$$
It is clear that $H^2$ is a $H^\infty$-module, with the action $f\cdot h$ of $f\in H^\infty$ on $h\in H^2$ given by pointwise multiplication of $f$ and $h$. 

In Section~\ref{sectionH2}, we show that 
 $H^2$ is a flat, but not faithfully flat $H^\infty$-module (Theorem~\ref{sarason}). 
 \label{pageref_sarason}It is also observed that for any finitely generated, nonzero, closed, proper ideal $\mathfrak{n}$ of $H^\infty$, we have $\mathfrak{n} H^2\subsetneq H^2$.
 
\subsection{Preliminaries}

 We recall some preliminaries that will be needed in both of the remaining sections. 
 
  Let $A$ be a commutative unital complex semisimple Banach algebra. The dual space $A^*$ of $A$ consists of all continuous linear complex-valued maps defined on $A$. The {\em maximal ideal
space} $M(A)$ of $A$  is the set of all nonzero homomorphisms $\varphi\!:\! A \!\to\!  {\mathbb{C}}$. As $M(A)$ is a subset of $A^*$, it inherits the  weak-$\ast$ topology of $A^*$, called the {\em Gelfand topology} on $M(A)$. 
 There is a bijective correspondence between $M(A)$ and the collection of all maximal ideals of $A$, namely,  each nonzero complex homomorphism $\varphi \in M(A)$ corresponds to a maximal ideal $\mathfrak{m}:=\ker \varphi $ of $A$.  The topological space $M(A)$ is a compact Hausdorff space. The set $M(A)$ is contained in the unit sphere of  the Banach space ${\mathcal{L}}(A,{\mathbb{C}})$ of all complex-valued continuous linear maps on $A$  with the operator norm, which is given by $\|\varphi\|=\sup_{a\in A, \;\|a\|\le 1} |\varphi(a)|$ for all $\varphi \in {\mathcal{L}}(A,{\mathbb{C}})$.  
 Let $C(M(A))$ denote the Banach algebra  of 
complex-valued continuous functions on $M(A)$ with pointwise operations and the supremum norm. 
 The {\em Gelfand transform} $\widehat{a}\in C(M(A))$ of an element $a\in A$ is  
 defined by $\widehat{a}(\varphi):=\varphi(a)$ for all $\varphi \in M(A)$.  
 
The following result plays a key role in the rest of the article. In particular, the sequence of growing weights $1/\sqrt{r_{n-1}}$ (see below) will be needed to construct an element of the maximal ideal $\mathfrak{m}$ of $R$ in our results about showing $\mathfrak{m} M=M$ for appropriate $\mathfrak{m}$ and $R$-modules $M$. 
 
\begin{proposition}
\label{olympiad}
Let $(a_n)_{n\in \mathbb{N}}$ be a sequence in ${\mathbb{C}}$ not having a finite support$,$ and  such that 
$
\textstyle \sum\limits_{n=1}^\infty |a_n|^2\!<\!\infty.
$ 
 Define\footnote{As $(a_n)_{n\in \mathbb{N}}$ does not have finite support, each $r_n>0$.}  $
\textstyle r_n\!=\!\sum\limits_{k=n+1}^\infty |a_k|^2$ $(>0)$ for all $n\in {\mathbb{N}}\cup\{0\}$.  
 Then $\sum\limits_{n=2}^\infty |a_n|^2\frac{1}{\sqrt{r_{n-1}}}\!<\!\infty$.
\end{proposition}

\smallskip 

\noindent {\em Proof}. We have $(r_n)_{n\in {\mathbb{N}}}$ is a decreasing sequence of positive reals, and
 $$
 \textstyle 
 \begin{array}{rcl}
\frac{|a_n|^2}{\sqrt{r_{n-1}}} =\frac{r_{n-1}-r_n}{\sqrt{r_{n-1}}}
=\frac{(\sqrt{r_{n-1}}-\sqrt{r_n})(\sqrt{r_{n-1}}+\sqrt{r_n})}{\sqrt{r_{n-1}}}
  &\le&  \frac{(\sqrt{r_{n-1}}-\sqrt{r_n})(\sqrt{r_{n-1}}+\sqrt{r_{n-1}})}{\sqrt{r_{n-1}}}\\
  &\le&  
2(\sqrt{r_{n-1}}-\sqrt{r_n}).
\end{array}
$$
 Thus the partial sums of $\sum\limits_{n=2}^\infty |a_n|^2\frac{1}{\sqrt{r_{n-1}}}$ form a Cauchy sequence:
 $$
 \textstyle 
\begin{array}{rcl}
 \sum\limits_{k=m+1}^n \frac{|a_k|^2}{\sqrt{r_{k-1}}}
  &\le & 
 2(\sqrt{r_{m}}\!-\!\sqrt{r_{m+1}}+\sqrt{r_{m+1}}\!-\!\sqrt{r_{m+2}}+\cdots+\sqrt{r_{n-1}}\!-\!\sqrt{r_n})
 \\
  &=& 2(\sqrt{r_m}-\sqrt{r_n})<2 \sqrt{r_m} <\epsilon,
 \end{array}
$$
 whenever $n>m>N$, where $N$ is such that $\sum\limits_{k=N+1}^\infty |a_k|^2<\frac{\epsilon^2}{4}$. 
 $\hfill\Box$

\noindent The organisation of the paper is as follows. In the next section, we prove our first set of main results (all about the $L^\infty(X,\mu)$-module $L^2(X,\mu)$), and in the last section, we establish the second set of main results (namely that the $H^\infty$-module $H^2$ is not faithfully flat). 

\vspace{-0.42cm}

\section{Unfaithfulness of the flat $L^\infty(X,\mu)$-module $L^2(X,\mu)$.}
\label{section_ell2}


\subsection{$L^2(X,\mu)$ is a flat $L^\infty(X,\mu)$-module} 


\begin{proposition}
\label{18_9_24_1810}
Let $X$ be a locally compact Hausdorff topological space and $\mu$ be a positive Radon measure on $X$. 
Then $L^2(X,\mu) $ is a flat $L^\infty(X,\mu)$-module.
\end{proposition}

\vspace{-0.1cm}

\noindent {\em Proof}. 
 Let $r_1 m_1+\cdots+r_nm_n=0\in L^2(X,\mu)$ for some $n\in {\mathbb{N}}$, and functions $r_1,\cdots, r_n\in L^\infty(X,\mu) $ and $m_1,\cdots, m_n \in L^2(X,\mu)$. Then there exists a subset $N$ of $X$ with $\mu(N)=0$ such that for all $x\in X\setminus N$, 
 $
 \textstyle 
r_1(x)m_1(x)+ \cdots + r_n(x)m_n(x)\!=\!0,
$  
 i.e., ${\mathbf{m}}(x):=(m_1(x),\cdots, m_n(x))\in {\mathbb{C}}^n$ is in the orthogonal complement ${\mathbf{r}}(x)^\perp$ of ${\mathbf{r}}(x):=(r_1(x), \cdots, r_n(x))\in {\mathbb{C}}^n$ with respect to the usual Euclidean inner product $\langle\cdot,\cdot\rangle$ (and corresponding norm denoted by $\|\cdot\|$ below) on ${\mathbb{C}}^n$. If ${\mathbf{r}}(x)\!\neq\! 0$, there exist $n-1$ orthonormal vectors ${\mathbf{e}}_1(x),\cdots, {\mathbf{e}}_{n-1}(x)$ forming a basis for ${\mathbf{r}}(x)^\perp$, and we set ${\mathbf{e}}_n(x)\!=\!{\mathbf{0}}\!\in\! {\mathbb{C}}^n$. If ${\mathbf{r}}(x)\!=\!0$, then let $\{{\mathbf{e}}_1(x),\cdots, {\mathbf{e}}_{n}(x)\}$ be any orthonormal basis for ${\mathbb{C}}^n$.  
 For $x\in X\setminus N$, define 
  $\mu_1(x),\cdots, \mu_{n}(x)$ by 
$\mu_j(x)=\langle {\mathbf{m}}(x), {\mathbf{e}}_j(x)\rangle$, $1\le j\le n$.  
For $x\in X\setminus N$, and $1\le i,j\le n$, set $\rho_{ij}(x)=\langle {\mathbf{e}}_j(x), e_i\rangle$, where $e_1,\cdots, e_n$ are the standard basis vectors for ${\mathbb{C}}^n$. Then by the Cauchy-Schwarz inequality, $|\rho_{ij}(x)|\le \|{\mathbf{e}}_j(x)\| \| e_i\|= 1$. Thus we have defined an element $\rho_{ij}$ of $ L^\infty(X,\mu) $ for each $1\le i,j,\le n$. Moreover, for all $x\in X\setminus N$, we have for $1\le j\le n$ that 
$$
\textstyle 
{\scaleobj{0.9}{
\begin{array}{rcl}
 0=\langle {\mathbf{r}}(x), {\mathbf{e}}_j(x)\rangle
  &=&  \left[\begin{smallmatrix} r_1(x) &\cdots & r_n(x)\end{smallmatrix}\right]
 \left[\begin{smallmatrix} \langle {\mathbf{e}}_j(x),e_1\rangle \\ \vdots \\ \langle {\mathbf{e}}_j(x),e_n\rangle \end{smallmatrix}\right]\\[0.42cm]
  &=& r_1(x) \rho_{1j}(x)+\cdots +r_{n}(x) \rho_{nj}(x).
 \end{array}}}
  $$
  Thus in $L^\infty(X,\mu)$, we have $r_1 \rho_{1j}+\cdots+r_n\rho_{nj}=0$, $1\le j\le n$. Moreover, $\mu_1,\cdots, \mu_n\in L^2(X,
  \mu) $ because for $1\le i\le n$, we have  
$$
\textstyle 
{\scaleobj{0.9}{
\begin{array}{rcl}
\int_{X\setminus N} |\mu_i(x)|^2 \;d\mu(x)
 &=& 
\int_{X\setminus N} |\langle {\mathbf{m}}(x), {\mathbf{e}}_i(x)\rangle|^2 \;d\mu(x)\\[0.1cm]
 &\le& 
\int_{X\setminus N} \|{\mathbf{m}}(x)\|^2 1 \;d\mu(x)
 =\|m_1\|_2^2+\cdots+\|m_n\|_2^2.
\end{array}}}
$$
Finally,  for $x\in X\setminus N$, we have 
$$
\textstyle 
{\scaleobj{0.9}{
\begin{array}{rcl}
\sum\limits_{i=1}^{n} m_i(x) e_i
 &=& {\mathbf{m}}(x)
=\sum\limits_{j=1}^{n} \langle {\mathbf{m}}(x),{\mathbf{e}}_j(x)\rangle\;\! {\mathbf{e}}_j(x) 
= \sum\limits_{j=1}^{n}  \mu_j(x) {\mathbf{e}}_j(x) \\
 &=& 
\sum\limits_{j=1}^{n}  \mu_j(x) 
\sum\limits_{i=1}^{n} \langle {\mathbf{e}}_j(x) , e_i\rangle e_i
= \sum\limits_{i=1}^{n} \sum\limits_{j=1}^{n}  \rho_{ij}(x) \mu_j(x)e_i.
\end{array} }}
$$
and so,  we get $
\textstyle 
m_i =\sum\limits_{j=1}^{n} \rho_{ij} \mu_j$ in $L^2(X,\mu)$ for all $1\le i\le n$.
$\hfill\Box$

\subsection{The case of finitely generated maximal ideals}

We claim that  for any finitely generated, nonzero, proper ideal $\mathfrak{n}$ of $L^\infty(X,\mu)$, we have $\mathfrak{n} L^2(X,\mu)\subsetneq L^2(X,\mu)$. To do this, we first show that $L^\infty(X,\mu)$ is a B\'ezout ring. Recall that a commutative ring is called {\em B\'ezout} if every finitely
generated ideal is principal.

\begin{lemma}
 \label{prop_Bezout}
 Every finitely generated ideal in $L^\infty(X,\mu)$ is principal$,$ that
 is$,$ $L^\infty(X,\mu)$ is B\'ezout ring.
\end{lemma}

\noindent Before we give the proof, we collect some useful
observations first. For 
$f\in L^\infty(X,\mu)$, let $
|f|, \overline{f}\in L^\infty(X,\mu)$ be the functions obtained by taking pointwise complex absolute value, and pointwise complex conjugation, respectively. 
Then for a representative function $f:X\to {\mathbb{C}}$ of an element in $L^\infty(X,\mu)$,   $f=|f|\cdot u_f$, where $u_f\in L^\infty(X,\mu)$ is given by 
$$
\textstyle 
u_f(x)=\left\{\begin{array}{cl}
            \frac{f(x)}{|f(x)|} & \textrm{if }\; f(x)\neq 0,\\[0.21cm]
            1 & \textrm{if }\;f(x)=0. 
          \end{array}\right. 
$$
We have $f\in L^\infty(X,\mu)$ if and only if
$|f|\in L^\infty(X,\mu)$.  
Also, $f\in L^\infty(X,\mu)$ if and only if
$\overline{f}\in  L^\infty(X,\mu)$. Also,
$u_f u_{\overline{f}}={\mathbf{1}}$ (the constant function, taking
value $1$ everywhere on $X$) and $|f|=f \overline{u_f}$.

In a commutative ring $R$, the ideal in $R$ generated by $r_1,\cdots, r_n$, $n\in {\mathbb{N}}$,  is denoted by $\langle r_1,\cdots, r_n\rangle$.
 
 \smallskip 
 
\noindent {\em Proof}. 
  It suffices to show that an ideal $\langle f,g\rangle$
  generated by $f,g\!\in\! L^\infty(X,\mu)$ is principal.  We claim 
  that $\langle f,g\rangle =\langle |f|+|g|\rangle.  $ 
  Since $\overline{u_f}, \overline{u_g}\in L^\infty(X,\mu)$, we have
  $|f|+|g|=f \overline{u_f}+g \overline{u_g}\in \langle
  f,g\rangle.  $
  Thus $\langle |f|+|g|\rangle\subset \langle f,g\rangle $.
 To show the reverse inclusion, let us  define $F$ by
$$
\textstyle 
F(x) =\left\{\begin{array}{cl}
    \frac{f(x)}{|f(x)|+|g(x)|} & 
    \textrm{if }\; |f(x)|\!+\!|g(x)|\!\neq \!0,\\[0.21cm]

    1 & \textrm{if }\; |f(x)|\!+\!|g(x)|\!=\!0, 
           \end{array}\right. 
$$
$x\in X$. Then $|F(x)|\leq 1$ for all $x\in X$, and so
$F \in L^\infty(X,\mu)$. Moreover, $f=F \cdot (|f|+|g|)$,
and so $f\in \langle |f|+|g|\rangle$.  Similarly, $g\in
\langle |f|+|g|\rangle$ too. Hence $\langle f,g\rangle
\subset \langle |f|+|g|\rangle$.
 Consequently, $\langle f,g\rangle =\langle
|f|+|g|\rangle$. 
$\hfill\Box$ 

\smallskip

\noindent 
If  the nonzero maximal ideal $\mathfrak{m}$ of $L^\infty(X,\mu)$ is finitely generated, then is $\mathfrak{m}L^2(X,\mu)\subsetneq L^2(X,\mu)$? We claim that  for any finitely generated, nonzero, proper ideal $\mathfrak{n}$ of $L^\infty(X,\mu)$, $\mathfrak{n} L^2(X,\mu)\subsetneq L^2(X,\mu)$. 

\begin{proposition} 
\label{19_9_24_1616}
If $\mu$ is a $\sigma$-finite positive Radon measure on a locally compact Hausdorff space $X,$ 
then for every nonzero$,$ finitely generated$,$ proper ideal $\mathfrak{n}$ of $L^\infty(X,\mu),$  
 we have  $\mathfrak{n}L^2(X,\mu)\subsetneq L^2(X,\mu)$. 
\end{proposition}

\noindent {\em Proof}. As $L^\infty(X,\mu)$ is a B\'ezout ring, $\mathfrak{n}=\langle g\rangle$ for some $g\in L^\infty(X,\mu)$, and as $\mathfrak{n}$ is nonzero, so is $g$. Then 
$$
\textstyle 
\mathfrak{n}L^2(X,\mu)=\langle g\rangle L^2(X,\mu) 
=gL^\infty(X,\mu) L^2(X,\mu) = gL^2(X,\mu).
$$  
We claim that $gL^2(X,\mu)\subsetneq L^2(X,\mu)$. 
Suppose on the contrary that $gL^2(X,\mu)=L^2(X,\mu)$. 

We will first show that $gL^2(X,\mu)=L^2(X,\mu)$ implies that $g\neq 0$ almost everywhere in $X$. Indeed, if $\mu\{x\in X: g(x)=0\}>0$, then since $\mu$ is a Radon measure, there exists a compact set  $K\subset \{x\in X\!:\! g(x)\!=\!0\}$ with $0<\mu(K)<\infty$.  Let $f=\mathbf{1}_K$ be the indicator function of $K$ (i.e., pointwise equal to $1$ on $K$ and $0$ on $X\setminus K$). Then clearly $f\in L^2(X,\mu)$, and so (thanks to our assumption that $gL^2(X,\mu)=L^2(X,\mu)$) there exists an $h\in L^2(X,\mu)$ such that $f=gh$. But at any point  $x\in K$, we then get $1=\mathbf{1}_K(x)=f(x)=g(x)h(x)=0h(x)=0$, a contradiction. 

As $g(x)\neq 0$ almost everywhere in $X$, it follows that the multiplication map $M_g: L^2(X,\mu) \to L^2(X,\mu)$, $h\mapsto gh$,  is injective. But thanks to our assumption $gL^2(X,\mu)=L^2(X,\mu)$, we also know that $M_g$ is surjective. Thus $M_g$ is invertible. As $\mu$ is $\sigma$-finite, we conclude that  $g$ is invertible as an element of $L^\infty(X,\mu)$, that is, $g^{-1}\in L^\infty(X,\mu)$ (see, e.g., \cite[Problem~67]{Hal}). But then $\mathfrak{n}=\langle g\rangle=L^\infty(X,\mu)$, contradicting the properness of $\mathfrak{n}$. 
$\hfill\Box$

\subsection*{Ultrafilters and $M(\ell^\infty)$}

\noindent In order to show the non-faithfulness, we will use the maximal ideal space of $\ell^\infty$, which is known to be in correspondence with the set of ultrafilters on ${\mathbb{N}}$. We begin by recalling some preliminaries about ultrafilters. 

Let $\ell^\infty$ be the Banach algebra of all complex bounded sequences 
on ${\mathbb{N}}=\{1,2,3,\cdots\}$ with pointwise operations and the supremum norm, which is given by 
$$
\textstyle
\|f\|_\infty:=\sup_{n\in {\mathbb{N}}} |f(n)|\; \text{ for }\;f=(f(n))_{n\in {\mathbb{N}}} \in \ell^\infty.
$$
Let $\ell^2 $ be the Hilbert space of all square summable complex sequences with the inner product corresponding to the norm $\|\cdot\|_2$, where 
$$
\textstyle
\|h\|_2^2=\sum\limits_{n=1}^\infty |h(n)|^2 \;\text{ for all }\; h=(h(n))_{n\in {\mathbb{N}}} \in \ell^2 .
$$
Then with the action of $\ell^\infty $ on $\ell^2 $ given by termwise multiplication, 
 $
\textstyle 
f\cdot h=\pmb{(}\;\!f(n)h(n)\pmb{)}_{n\in {\mathbb{N}}}$  for $ f\in \ell^\infty$, $h\in \ell^2$,  $\ell^2 $ is an $\ell^\infty $-module.  As a consequence of our results in this section, it will follow that as an $\ell^\infty$-module, $\ell^2$ is  flat,  but not faithfully flat (see Example~\ref{example_ellinfty}). 

For $n\in {\mathbb{N}}$, define $\varphi_n: \ell^\infty \to {\mathbb{C}}$ by $\varphi_n(f)=f(n)$ for all $f\in \ell^\infty $. 
Point evaluations $\varphi_n$, $n\in {\mathbb{N}}$, are nonzero complex homomorphisms. Thus ${\mathbb{N}}$ can be identified  
 with a subset of $M(\ell^\infty )$ (and the topology on ${\mathbb{N}}$, induced from the Gelfand topology on $M(\ell^\infty )$, coincides with the usual Euclidean topology on ${\mathbb{N}}$ as a subset of ${\mathbb{R}}$).  For each $n\in {\mathbb{N}}$, the maximal ideal $\mathfrak{m}_n:=\ker \varphi_n$ has  the property that $\mathfrak{m}_n\ell^2\subsetneq \ell^2$. Indeed, for example, if $f\in \ell^2$ is the sequence with $1$ as the $n$th term, and all other terms $0$, then $f\not\in  \mathfrak{m}_n\ell^2$ (otherwise $f\!=\!gh$, for some $g\!\in\! \mathfrak{m}_n$, $h\!\in\! \ell^2$,   and so  $0=\widehat{g}(\varphi_n)=g(n)$, giving $1=f(n)=g(n)h(n)=0h(n)=0$,  a contradiction). Despite the fact that $\mathfrak{m}_n\ell^2\subsetneq \ell^2$ for each $n\in {\mathbb{N}}$, there exist many maximal ideals $\mathfrak{m}$ of $\ell^\infty$ for which $\mathfrak{m}\ell^2=\ell^2$. 
 We recall that 
  $M(\ell^\infty )$ can be identified with the set of all ultrafilters on ${\mathbb{N}}$ as described below. 
  We refer to the article \cite{KomTot} for background on ultrafilters. 
  
  An {\em ultrafilter ${\mathcal{F}}$ on ${\mathbb{N}}$} is a collection of subsets of ${\mathbb{N}}$ with the properties listed below:
  
  (F1) $\;$ ${\mathbb{N}}\in {\mathcal{F}}$, $\emptyset \not \in {\mathcal{F}}$. 
  
  (F2) $\;$ If $A\in {\mathcal{F}}$ and $A\subset  B$, then $B\in {\mathcal{F}}$. 
  
  (F3) $\;$ If $A,B\in {\mathcal{F}}$, then $A\cap B\in {\mathcal{F}}$. 
 
  (F4) $\;$ $A\in {\mathcal{F}}$ if and only if ${\mathbb{N}}\setminus A \not\in{\mathcal{F}}$. 
  
  \noindent 
  A family of subsets of ${\mathbb{N}}$ with the properties (F1), (F2), (F3) is called a {\em filter on ${\mathbb{N}}$}. 
  For a fixed $n\in {\mathbb{N}}$, the ultrafilter 
  $$
    {\mathcal{F}}_n=\{A\subset  {\mathbb{N}}: n\in A\} 
    $$ 
    is called a {\em principal ultrafilter on ${\mathbb{N}}$}. 
    As ${\mathbb{N}}$ is an infinite set, it can be shown (using Zorn's lemma) that there exists a non-principal ultrafilter on ${\mathbb{N}}$. 
  If ${\mathcal{F}}$ is a non-principal ultrafilter on ${\mathbb{N}}$, then it contains all cofinite sets in ${\mathbb{N}}$ (i.e., all subsets of ${\mathbb{N}}$ whose complements are finite sets). For a bounded complex sequence $(a_n)_{n\in {\mathbb{N}}}$, a complex number $L$,  and an ultrafilter ${\mathcal{F}}$ on ${\mathbb{N}}$,   we write 
  $$
  \textstyle 
  \lim\limits_{\mathcal{F}} a_n=L
  $$
  if for every $\epsilon\!>\!0$, the set $\{n\in {\mathbb{N}}: |a_n-L|<\epsilon\}$ belongs to ${\mathcal{F}}$. It is clear that if a complex sequence $
  (a_n)_{n\in {\mathbb{N}}}$ is convergent in ${\mathbb{C}}$ with limit $L$, then for each non-principal ultrafilter ${\mathcal{F}}$ on ${\mathbb{N}}$, we have 
  $$
  \textstyle \lim\limits_{\mathcal{F}} a_n\!=\!L
  $$
   as well. For $m\in {\mathbb{N}}$, and the principal ultrafilter ${\mathcal{F}}_m$ on ${\mathbb{N}}$, we have 
   $$
   \textstyle \lim\limits_{{\mathcal{F}}_m}a_n =L
   $$
    if and only if $a_m=L$. 
  The maximal ideal $\mathfrak{m}_{\mathcal{F}}$ of $\ell^\infty$ corresponding to an ultrafilter ${\mathcal{F}}$ on ${\mathbb{N}}$ is 
  $$
  \textstyle
  \mathfrak{m}_{\mathcal{F}}=\{f\in \ell^\infty : \lim\limits_{\mathcal{F}} f(n)=0\}.
  $$
  This exhausts $M(\ell^\infty)$. See, e.g., \cite[Thm, \S2.5]{GilJer} (where only the real-valued case is considered, but the proof for the complex-valued case is the same, mutatis mutandis). 
  From the above, we know that if $m\in {\mathbb{N}}$, then $\mathfrak{m}_{{\mathcal{F}}_m} \ell^2  \subsetneq \ell^2 $. 
 With these preliminaries 
 about $M(\ell^\infty)$ in hand, we show the non-faithfulness of the flat $L^\infty(X,\mu)$-module $L^2(X,\mu)$ in the following subsection.

\subsection{$L^2(X,\mu)$ is not a faithfully flat $L^\infty(X,\mu)$-module}

Finally, we show the main result of this section. We will make a mild assumption, 
which will be satisfied in all the examples of our interest.

\begin{theorem}
\label{19_9_24_1617}
Let $X$ be a locally compact Hausdorff topological space$,$ $\mu$ be a positive Radon measure on $X,$ such that  there is a family of Borel sets $U_n$ of $X,$ $n\in {\mathbb{N}},$  such that 
\begin{itemize}
\item $U_1\subset U_2\subset U_3\subset \cdots,$ and $X=\textstyle \bigcup\limits_{n\in {\mathbb{N}}} U_n,$  
\item for all $n\!\in\! {\mathbb{N}},$ $\overline{U_n}$ is compact$,$ 
and   $\mu(U_{n}\setminus U_{n-1})>0,$ where $U_0:=\emptyset$. 
\end{itemize}
Let ${\mathcal{F}}$ be any non-principal ultrafilter on ${\mathbb{N}},$ and for each $n\in {\mathbb{N}},$ let $\varphi_n\in M\pmb{(}L^\infty(U_n\setminus U_{n-1}, \mu)\pmb{)}$.  
Define $\varphi\in M\pmb{(} L^\infty(X,\mu)\pmb{)}$ by 
$$
\textstyle 
\varphi(f)=\lim\limits_{\mathcal{F}} \varphi_n(f|_{U_n\setminus U_{n-1}}) \;\text{ for all }\; f\in L^\infty(X,\mu).
$$
If $\mathfrak{m}:=\ker \varphi,$ then $\mathfrak{m}L^2(X,\mu)=L^2(X,\mu)$. 
\end{theorem}

\noindent Note that under the above assumptions, $X$ is necessarily an infinite set. As each $U_n\setminus U_{n-1}$ is a Borel set, the restriction of the Radon measure $\mu$ to $U_n\setminus U_{n-1}$ is again a Radon measure (see, e.g., \cite[Chap.7, Ex.7]{Fol}). For each $n\in {\mathbb{N}}$, $L^\infty(U_n\setminus U_{n-1}, \mu)$ is a Banach algebra with pointwise operations and the (essential) supremum norm. If the complex homomorphisms $\varphi_n\in M\pmb{(}L^\infty(U_n\setminus U_{n-1}, \mu)\pmb{)}$, $n\in {\mathbb{N}}$, then for any $g\in L^\infty(X,\mu)$, 
the complex sequence $\pmb{(}\varphi_n(g|_{U_n\setminus U_{n-1}})\pmb{)}_{n\in {\mathbb{N}}}\in \ell^\infty$ because for all $n\in {\mathbb{N}}$, 
 $
\noindent 
|\varphi_n(g|_{U_n\setminus U_{n-1}})|
\le \|\varphi_n\|\|g|_{U_n\setminus U_{n-1}}\|_\infty 
\le 1\|g\|_\infty.
$ 
Here $\|\varphi_n\|$ denotes the operator norm of $\varphi_n:L^\infty(U_n\setminus U_{n-1}, \mu)\to {\mathbb{C}}$ 
which is equal to $1$ (see, e.g., \cite[Prop.~2.22]{Dou}).  
If ${\mathcal{F}}$ is any non-principal ultrafilter on ${\mathbb{N}}$ and $\varphi_n\in M\pmb{(}L^\infty(U_n\setminus U_{n-1}, \mu)\pmb{)}$, $n\in {\mathbb{N}}$, then 
 $\varphi: L^\infty(X,\mu)\to {\mathbb{C}}$, as defined in the theorem statement above, is a continuous linear transformation, and it respects multiplication. If $\mathbf{1}$ denotes the constant function on $X$ taking value $1$ everywhere, then $\varphi(\mathbf{1})=1$. So $\varphi \in M\pmb{(}L^\infty(X,\mu)\pmb{)}$. 

\medskip 

\noindent {\em Proof}.  Let $f\in L^2(X,\mu)$. We will construct a function $w:X\to (0,+\infty)$  such that 
$w(x)\to +\infty$ as $x\to \boldsymbol{\infty}$, and $fw\in \ell^2 $. Here $X\cup \{\boldsymbol{\infty}\}$ is the one-point/Alexandroff  compactification of $X$. Below, by $f\in L^2(X,\mu)$ having a {\em compact support}, we mean that there exists a representative of $f$ that is identically $0$ outside a compact set. 
\begin{itemize}
\item[$1^\circ$] If $f$ has compact support, then we simply set $w|_{U_1}\!=\!1$ and 
 for $n\ge 2$, $w|_{U_n\setminus U_{n-1}}=n$.  

\item[$2^\circ$] Suppose $f$ does not have compact support. We set $w|_{U_1}\!=\!1$ and 
 for all $n\ge 2$, $w|_{U_n\setminus U_{n-1}}=\frac{1}{\sqrt[4]{r_{n-1}}}$, 
 where for $n\in {\mathbb{N}}$,
  $$
  \textstyle
 r_n\!:=\!\sum\limits_{k=n+1}^\infty a_k^2, \;\text{ and }\;
 a_k^2\!:=\!\int_{x\in U_k\setminus U_{k-1}} |f(x)|^2 d\mu(x) \;\text{ for }\;k\in {\mathbb{N}}.
$$
Note that $r_n>0$ since  $f$ does not have compact support. We have $r_n\to 0$ as $n\to +\infty$ because $f\in L^2(X,\mu)$.  
Thus $w(x)\to +\infty$ as $x\to {\boldsymbol{\infty}}$. Moreover, using Proposition~\ref{olympiad}, we now show that $fw\in L^2 (X,\mu)$. We have 
$$
\begin{array}{rcl}
 \quad \quad\quad&&  \int_X |f(x)|^2 w(x)^2d\mu(x)\\[0.18cm]
 &=& 
\int_{U_1} |f(x)|^2d\mu(x)+
\sum\limits_{n=2}^\infty 
\int_{x\in U_{n}\setminus U_{n-1}} |f(x)|^2 \frac{1}{\sqrt{r_{n-1}}} d\mu(x)     
  \quad \quad\quad(\star)\\
 &\le & \|f\|_2^2 + \sum\limits_{n=2}^\infty \frac{a_n^2}{\sqrt{r_{n-1}}}<\infty.
 \end{array}
$$
\end{itemize}
Define $g\in L^\infty(X,\mu)$ by 
$$
\textstyle 
g(x)=\frac{1}{w(x)}\;\text{ for all }\;x\in X.
$$
 Then $g(x)\to 0$ as $x\to \boldsymbol{\infty}$. Thus as $w|_{U_n\setminus U_{n-1}}$ is a constant function, taking value which we denote by $\omega_n$, and since by construction,  $(\omega_n)_{n\in {\mathbb{N}}}$ diverges to $+\infty$ as $n\to +\infty$, we get 
 $$
 \textstyle 
 \varphi(g)=\lim\limits_{\mathcal{F}} \varphi_n(g|_{U_n\setminus U_{n-1}})=
 \lim\limits_{\mathcal{F}} \frac{1}{\omega_n}=0.
 $$
 So $g\in \ker \varphi=\mathfrak{m}$. Next, define $h$  by 
 $$
 \textstyle h(x)=\frac{f(x)}{g(x)}\;\text{ for all }\;x\in X.
 $$
  By ($\star$) above,  $h$ belongs to $L^2(X,\mu)$. So $f=gh$, where $g\in \mathfrak{m}$ and $h\in L^2(X,\mu) $. 
  This completes the proof. 
$\hfill\Box$ 

\medskip

\noindent An immediate consequence of Theorem~\ref{19_9_24_1617} is that, 
under the same assumptions on $X$ and $\mu$, 
the $L^\infty(X,\mu)$-module $L^2(X,\mu)$ is not faithfully flat. 

\subsection{Applications/examples}

The following list of examples is motivated by the consideration of signal spaces arising in control theory, see, e.g., \cite{Qua}. 

\begin{example}
\label{example_ellinfty}
\em{ 
Let $X={\mathbb{N}}$ be endowed with the usual Euclidean topology induced from ${\mathbb{R}}$. Then the compact subsets of ${\mathbb{N}}$ are finite.  Let $\mu$ be the counting Radon measure, defined by setting for any compact subset $K$ of ${\mathbb{N}}$, $\mu(K)$ to be the number of elements in $K$. Then $L^\infty(X,\mu)=\ell^\infty$ and $L^2(X,\mu)=\ell^2$. 
 Thus by Proposition~\ref{18_9_24_1810}, $\ell^2$ is a flat $\ell^\infty$-module. 
 As $\mu$ is $\sigma$-finite, it follows from Proposition~\ref{19_9_24_1616} that 
 for every finitely generated, nonzero, proper ideal $\mathfrak{n}$ of $\ell^\infty$, we have $\mathfrak{n}\ell^2 \subsetneq \ell^2$. Finally, taking $U_n=\{1,2,\cdots, n\}$, we have
 \begin{itemize}
\item $U_1\subset U_2\subset U_3\subset \cdots$, and ${\mathbb{N}}= \bigcup\limits_{n\in {\mathbb{N}}} U_n$, 

\item for all $n\in {\mathbb{N}}$, $\overline{U_n}$ is compact, 
and   $\mu(U_{n}\setminus U_{n-1})=\mu\{n\}=1>0$. 
\end{itemize}
Thus $\ell^2$ is not a faithfully flat $\ell^\infty$-module by Theorem~\ref{19_9_24_1617}.
}
\hfill$\Diamond$
\end{example}

\begin{example}
\em{
 Let $X={\mathbb{R}}$ be endowed with the usual Euclidean topology, 
 and $\mu$ be the Lebesgue measure $m$.  Then $L^2({\mathbb{R}},m)$ is a flat $L^\infty({\mathbb{R}},m)$-module. 
 As $m$ is $\sigma$-finite, 
 it follows that for every finitely generated, nonzero, proper ideal $\mathfrak{n}$ of $L^\infty({\mathbb{R}},m)$,  $\mathfrak{n}L^2({\mathbb{R}},m) \subsetneq L^2({\mathbb{R}},m)$. Finally, taking $U_n=(-n,n)$, we have 
 \begin{itemize}
\item $U_1\subset U_2\subset U_3\subset \cdots$, and ${\mathbb{R}}= 
 \bigcup\limits_{n\in {\mathbb{N}}}  U_n$,

\item for all $n\in {\mathbb{N}}$, $\overline{U_n}=[-n,n]$ is compact, 
and  $\mu(U_{n}\setminus U_{n-1})>0$. 
\end{itemize}
Thus $L^2({\mathbb{R}},m)$ is not a faithfully flat $L^2({\mathbb{R}},m)$-module. }
\hfill$\Diamond$
\end{example}

\begin{example}
\em{
 Let ${\mathbb{T}}=\{z\in {\mathbb{C}}:|z|=1\}$ and $X={\mathbb{T}}$ be endowed with the usual Euclidean topology induced from ${\mathbb{R}}^2$, 
 and $\mu$ be the Lebesgue measure $m$ on ${\mathbb{T}}$.  Then $L^2({\mathbb{T}},m)$ is a flat $L^\infty({\mathbb{T}},m)$-module. 
 As $m$ is $\sigma$-finite,  it follows that for every finitely generated, nonzero, proper ideal $\mathfrak{n}$ of $L^\infty({\mathbb{T}},m)$, we have $\mathfrak{n}L^2({\mathbb{T}},m) \subsetneq L^2({\mathbb{T}},m)$. Finally, taking $U_n=\{e^{i\theta}: \frac{1}{n}<|\theta|\le \pi\}\cup\{1\}$, we have 
 \begin{itemize}
\item $U_1\subset U_2\subset U_3\subset \cdots$, and ${\mathbb{T}}= \bigcup\limits_{n\in {\mathbb{N}}} U_n$,

\item  for all $n\in {\mathbb{N}}$, $\overline{U_n}$ is compact, 
and   $\mu(U_{n}\setminus U_{n-1})>0$. 
\end{itemize}
Thus $L^2({\mathbb{T}},m)$ is not a faithfully flat $L^2({\mathbb{T}},m)$-module. }
\hfill$\Diamond$
\end{example}

\section{$H^2$ is a flat, but not faithfully flat $H^\infty$-module.}
\label{sectionH2}
   
 \subsection{$H^2$ is a flat $H^\infty$-module} 
 The vector-valued Beurling theorem (the Lax-Halmos theorem, see, e.g., \cite[Cor. 6, pp. 17-18]{Nik}),  characterising shift-invariant subspaces of $H^2$, implies that $H^2$ is a flat $H^\infty$-module. This result was also shown in \cite[Prop. 8]{Qua} (in the half-plane setting), but we include a short proof for completeness. Below,  $(H^2)^n$ is the Hilbert space which is the direct sum of $n$ ($\in {\mathbb{N}}$) copies of $H^2$. 

\begin{proposition}
$H^2$ is a flat $H^\infty$-module. 
\end{proposition}

\noindent {\em Proof}. 
Let $n\in {\mathbb{N}}$ and $r_1,\cdots, r_n\in H^\infty$. Consider the  
${\mathbb{C}}$-linear  map ${\mathbf{r}}:(H^2)^n\to H^2$ given by ${\mathbf{r}}(h_1,\cdots, h_n) =r_1 h_1+\cdots+r_n h_n$ for all $ (h_1,\cdots, h_n)\in (H^2)^n$. It is clear that $\ker {\mathbf{r}}$ is a shift-invariant subspace (i.e., $(zh_1,\cdots, zh_n)\in \ker {\mathbf{r}}$ whenever $(h_1,\cdots, h_n)\in \ker {\mathbf{r}}$). By the Beurling-Lax-Halmos theorem, there exists a  $k\in {\mathbb{N}}$ and a matrix $[\rho_{ij}]$ of size $n\times k$ with $H^\infty$  entries, such that  
$$
\textstyle 
 \ker {\mathbf{r}}=\Big\{\sum\limits_{j=1}^k \rho_{ij} \varphi_j: \varphi_1, \cdots, \varphi_k\in H^2\Big\}.
 $$
 For a $j\!\in\! \{1,\cdots, k\}$, taking $\varphi_{i}\!=\!\delta_{ij}$ (equal to the constant function ${\mathbf{0}}$ if $i \neq j$ and the constant function ${\mathbf{1}}$ if $i\!=\!j$), $1\!\le\! i\!\le\! k$, we get $(\rho_{1j},\cdots, \rho_{nj})\!\in\! \ker {\mathbf{r}}$. So $H^2$ is a flat $H^\infty$-module.
 $\hfill\Box$

\subsection{The case of finitely generated maximal ideals}
\label{Mortini_subsec}

 We claim that  for any finitely generated, nonzero, closed, proper ideal $\mathfrak{n}$ of $H^\infty$, we have $\mathfrak{n} H^2\subsetneq H^2$. 
   Using a result from \cite{Mor} (saying that any  finitely generated, nonzero ideal of $H^\infty$ is closed if and only if  it is a principal ideal generated by  an inner function), $\mathfrak{n}=\langle g\rangle$, where $g\in H^\infty$ is an inner function. 
Hence $\mathfrak{n} H^2=\langle g\rangle H^2=gH^\infty H^2=  gH^2 \subsetneq H^2$. 
(That the last inclusion is strict: Otherwise, in particular ${\mathbf{1}}\in H^2=gH^2$, and so by looking at boundary values on ${\mathbb{T}}$, the pointwise complex conjugate $\overline{g}$ of $g$ satisfies  $\overline{g}=\frac{1}{g} \in H^2$. Thus $g,\frac{1}{g}\in H^2$, implying that the Fourier coefficients $\hat{g}_m$ of $g$ are zero for all $m\neq 0$, i.e., $g=\alpha$ with $|\alpha|=1$. But then $\mathfrak{n}=\langle g\rangle=H^\infty$, contradicting the properness of $\mathfrak{n}$.)

 \subsection{Some background on Hardy spaces} 
 
 We recall a few facts we need for  the proof of the non-faithfulness of the flat $H^\infty$-module $H^2$.  Background on Hardy spaces can be found, e.g., in \cite{Gar} and \cite{Hof}. 
   
 We recall that an {\em inner function} is a $g\in {\mathcal{O}}({\mathbb{D}})$ such that $|g(z)|\le 1$ for all $z\in {\mathbb{D}}$ and such that $|g(e^{i\theta})|=1$ for almost all $\theta \in (-\pi,\pi]$. An {\em outer function} is an analytic function $F\in {\mathcal{O}}({\mathbb{D}})$ having the form 
$$
\textstyle 
F(z)=\alpha \exp\big( \frac{1}{2\pi} \int_{-\pi}^\pi \frac{e^{i\theta}+z}{e^{i\theta}-z} k(\theta) d\theta\big)
$$
where $k:{\mathbb{T}}\to {\mathbb{R}}$ belongs to $L^1({\mathbb{T}})$ and $\alpha \in {\mathbb{T}}$. Then $k(\theta)=\log |F(e^{i\theta})|$ for almost all $\theta \in (-\pi,\pi]$. 

For the results claimed below, we refer to \cite[pp. 160-162]{Hof}. Consider the identity function ${\mathbf{z}} \in H^\infty$, given by ${\mathbb{D}}\owns z\mapsto z$. Then the map $\pi: M( H^\infty)\to {\mathbb{C}}$, 
 $
\pi(\varphi)=\varphi({\mathbf{z}})\text{ for all }\varphi \in  M( H^\infty),
$ 
 is a continuous map onto the closed  unit disk $\overline{{\mathbb{D}}}$  in ${\mathbb{C}}$. Over the open unit disk ${\mathbb{D}}$, $\pi$ is one-to-one, and maps ${\mathbb{D}}$ homeomorphically onto an open subset of $M( H^\infty)$, be sending $\lambda \in {\mathbb{D}}$ to $\varphi_\lambda \in M(H^\infty)$, where  $\varphi_\lambda$ is  given by $\varphi_\lambda (f)=f(\lambda)$ for all $f\in H^\infty$. 
   The remainder of $M( H^\infty)$ is mapped by $\pi$ onto the unit circle. If $|\alpha|=1$, then $\pi^{-1}\{\alpha\}$ is the {\em fibre of $M( H^\infty)$ over $\alpha$}. If $g\in  H^\infty$ and $\alpha \in {\mathbb{T}}$, then the range of $\widehat{g}$ on the fibre $\pi^{-1}\{\alpha\}$ is the set of all $\zeta \in {\mathbb{C}}$ such that there exists a sequence $(\lambda_n)_{n\in {\mathbb{N}}}$ in ${\mathbb{D}}$ with the properties that 
  $
\lim\limits_{n\to \infty}\lambda_n=\alpha \text{ and } 
\lim\limits_{n\to \infty} g(\lambda_n)=\zeta.
$

It is clear that for each point evaluation $\varphi_\lambda$, where $\lambda\in {\mathbb{D}}$, the corresponding maximal ideal $\mathfrak{m}_\lambda:=\ker \varphi_\lambda$ of $H^\infty$  has the property that $\mathfrak{m}_\lambda  H^2\subsetneq H^2$ (because any $g\in \mathfrak{m}_\lambda$ has a zero at $\lambda$, and so each element of $\mathfrak{m}_\lambda H^2$ will also have a zero at $\lambda$, but ${\mathbf{1}}\in H^2$ does not vanish at $\lambda$). 

However, $H^2$ is not faithfully flat, since we will show below that for any $\mathfrak{m}:=\ker \varphi$, where $\varphi\in M( H^\infty)\setminus \pi^{-1}({\mathbb{D}})$, we have $\mathfrak{m}  H^2=  H^2$. More explicitly, if $\alpha \in {\mathbb{T}}$, 
$\varphi \in \pi^{-1} \{\alpha\}$ and $\mathfrak{m}:=\ker \varphi$, then $\mathfrak{m}  H^2=  H^2$. We prove this just when $\alpha=1$, but rotational symmetry yields the result for arbitrary $\alpha \in {\mathbb{T}}$. 

\noindent The proof of Theorem~\ref{sarason} is based on the following email reply  (paraphrased on disc instead of the half-plane)  by Donald Sarason to the question (of faithfulness of the flat $H^\infty$-module $H^2$) put to him by Alban Quadrat, \cite{QuaSar}: 
\begin{quote}
Let $f\in  H^2$. Then $f=gh$, where $h\in  H^2$ and $g\in  H^\infty$, and $g(z)\to 0$ as $z\to 1$ (so $g$ belongs to every maximal ideal of $ H^\infty$ associated with a point in the maximal ideal space lying in the fibre above $1$).  To produce $g$ and $h$ take a positive function $w$ on ${\mathbb{T}}$ such that $w>1$, $wf\in L^2$ and $w(e^{i\theta})\to \infty$ as $\theta\to 0$. 
Let $g$ be the outer function whose modulus on ${\mathbb{T}}$ is $1/w$, and let  $h=f/g$. 
\end{quote}
\noindent We give a complete proof by making the outline above explicit below, 
since the construction of such a $w$ with the required properties, is not at all trivial. 

\begin{theorem}
\label{sarason}
\noindent If $\varphi \!\in\! \pi^{-1}\{1\},$ then for the maximal ideal $\mathfrak{m}\!=\!\ker \varphi$ of $  H^\infty,$ $\mathfrak{m}  H^2 \!=\! H^2$.
\end{theorem} 

\noindent {\em Proof}.  Let $f\!\in\!  H^2$ be nonzero.  Then by the F. and M. Riesz  theorem (see, e.g., \cite[Theorem~6.13]{Dou}), it follows that the boundary function $f\in L^2({\mathbb{T}})$ cannot be identically zero on any arc of ${\mathbb{T}}$ with a positive measure.  Without loss of generality, we may also assume that $\|f\|_2\le 1$. Let $a_n>0$ be defined by 
$$
\textstyle 
a_n^2=\int_{\frac{1}{n+1}\le |\theta| <\frac{1}{n}} |f(e^{i\theta})|^2 d\theta 
\text{ for all } n\in {\mathbb{N}}.
$$
Then $a_n\in (0,1)$, and also $(a_n)_{n\in {\mathbb{N}}}\in \ell^2$. For $n\in {\mathbb{N}}$, set 
$$
\textstyle 
r_n=\sum\limits_{k=n+1}^\infty a_k^2.
$$
 Then $r_n\in (0,1]$ (thanks to $\|f\|_2\le 1$ and as each $a_n>0$) for all $n\in {\mathbb{N}}$.  
Define $w$ as follows:
$$
\textstyle 
w(e^{i\theta})
=
\left\{
\begin{array}{ll}
1 & \text{if }\theta\in (-\pi,\pi]\setminus (-\frac{1}{2},\frac{1}{2}),\\[0.1cm]
\min\{\frac{1}{\sqrt[4]{r_{n-1}}},n\} & \text{if } \frac{1}{n+1}\le |\theta|<\frac{1}{n}, \;n\ge 2.
\end{array}\right.
$$
We note that $w$ is pointwise $\ge 1$ (as 
 $
\textstyle 
r_{n-1}\le \|f\|_2^2\le 1$). 
 Also,  we have that $w(e^{i\theta})\to +\infty$ as $\theta\to 0$. 

We wish to construct an outer function $g$ with modulus $\frac{1}{w}$ on ${\mathbb{T}}$, i.e., with $k:=\log |g|_{{\mathbb{T}}}|=\log \frac{1}{w}$. 
So we must check that the $w$ constructed above has the property that $k=\log \frac{1}{w}\in L^1({\mathbb{T}})$. We have 
$$
\textstyle 
\begin{array}{rcl}
\int_{-\pi}^\pi |\log \frac{1}{w(e^{i\theta})}| d\theta
 &=& 
\int_{\pi >|\theta|\ge \frac{1}{2}} |\log 1| d\theta + \sum\limits_{n=2}^\infty 
\int_{\frac{1}{n+1}\le |\theta| <\frac{1}{n}} \log w(e^{i\theta})  d\theta 
\\
 &\le &  0+ \sum\limits_{n=2}^\infty  \frac{2}{n(n+1)} \log n
\le 2 \sum\limits_{n=2}^\infty  \frac{\log n}{n^2} .
\end{array}
$$
But there exists an $N\in {\mathbb{N}}$ such that for all $n>N$, 
$\log n<\sqrt{n}$ because 
$$
\textstyle \lim\limits_{x\to \infty} \frac{\log x}{\sqrt{x}}=0.
$$
For $n>N$, $\frac{\log n}{n^2}<\frac{1}{n^{3/2}}$. As
$ \sum\limits_{n=1}^\infty \frac{1}{n^{3/2}}<\infty$, we get by comparison that 
$$
\textstyle  \sum\limits_{n=2}^\infty  \frac{\log n}{n^2}<\infty.
$$
  So from the above, we conclude that 
 $\log \frac{1}{w}\in L^1({\mathbb{T}})$. 

 Let $g$ be the outer function with modulus $\frac{1}{w}$ on ${\mathbb{T}}$, i.e., with $k=\log \frac{1}{w}$,
$$
\textstyle 
g(z)=\exp \big( \frac{1}{2\pi} \int_{-\pi}^\pi \frac{e^{i\theta}+z}{e^{i\theta}-z} k(\theta) \;\!d\theta\big).
$$
Then $g\in {\mathcal{O}}({\mathbb{D}})$, $\frac{1}{w}\in L^\infty({\mathbb{T}})$, and  $g\in H^\infty$.
 By the construction of $w$, we have $w(e^{i\theta})\to \infty$ as $\theta\to 0$, that is, 
 $g(e^{i\theta})=\frac{1}{w(e^{i\theta})}\to 0$ as $\theta\to 0$. By Lindel\"of's Theorem (see, e.g., \cite[Ex.~7(c),  pp.88-89]{Gar}), it follows that 
 $$
\quad \quad \quad \quad \quad \quad \quad \quad \quad \quad
\lim\limits_{{\mathbb{D}}\owns z\to 1} g(z)=0.
\quad \quad \quad\quad \quad \quad \quad \quad \quad \quad (\star\star)
 $$
 We claim that  $g\in \mathfrak{m}$, 
 that is $\varphi(g)=0$. This follows from the fact that the range of $\widehat{g}$ on $ \pi^{-1}\{1\}$ is 
 the set of all $\zeta \in {\mathbb{C}}$ such that there exists a sequence $(\lambda_n)_{n\in {\mathbb{N}}}$ in ${\mathbb{D}}$ satisfying 
$$
\lim\limits_{n\to \infty}\lambda_n=1 \; \text{ and } \;
\lim\limits_{n\to \infty} g(\lambda_n)=\zeta.
$$
But for each such sequence $(\lambda_n)_{n\in {\mathbb{N}}}$ in ${\mathbb{D}}$ converging to $1$, we have by ($\star\star$) that $(g(\lambda_n))_{n\in {\mathbb{N}}}$ converges to $0$. Hence the range of $\widehat{g}$ on $ \pi^{-1}\{1\}$ is just $\{0\}$. In particular $\widehat{g}(\varphi)=0$, that is, $\varphi(g)=0$. Hence $g\in \mathfrak{m}$. 

 As $g$ is given by an exponential, $g$ is never $0$, and so $\frac{1}{g}\in {\mathcal{O}}({\mathbb{D}})$. Define $h\in {\mathcal{O}}({\mathbb{D}})$ by $h=\frac{1}{g} f$. We claim that $h\in H^2$ by showing that  
 $wf\in  L^2({\mathbb{T}})$. We have 
 $$
 \textstyle 
\! \begin{array}{rcl}
 \int_{-\pi}^\pi \!|f(e^{i\theta})|^2 (w(e^{i\theta}))^2 d\theta 
  &\le & 
 \int_{\pi>|\theta|\ge \frac{1}{2}} \! |f(e^{i\theta})|^2d\theta 
 + \sum\limits_{n=2}^\infty 
 \int_{\frac{1}{n+1}\le |\theta| \le \frac{1}{n}}\! \frac{|f(e^{i\theta})|^2}{\sqrt{r_{n-1}}} d\theta 
 \\
  &\le&  \|f\|_2^2+ \sum\limits_{n=2}^\infty  \frac{a_n^2}{\sqrt{r_{n-1}}}<\infty,
 \end{array}
 $$
 where we used Proposition~\ref{olympiad} to obtain the last inequality. 
 Thus  $f=gh$ with $g\in \mathfrak{m}$ and $h\in H^2$, showing that $\mathfrak{m} H^2=H^2$. 
$\hfill\Box$

\subsection{The half plane case} 
In control theory, besides the Hardy spaces on the disc (corresponding to `discrete time systems' arising from difference equations), one also encounters Hardy spaces on the right half plane ${\mathbb{C}}_+=\{s\in {\mathbb{C}}: \text{Re}\;s>0\}$ (corresponding to `continuous time systems' arising from differential equations). We observe that the results above also hold in this case, as explained below. 

The {\em Hardy Hilbert space} $H^2({\mathbb{C}}_+)$ is the set of all  $F\in {\mathcal{O}}({\mathbb{C}}_+)$ with 
$$
\textstyle
\|F\|_2^2:=\sup\limits_{x>0} \int_{-\infty}^\infty |F(x+iy)|^2\;\! dy<\infty.
$$
Then $H^2({\mathbb{C}}_+)$ is a Hilbert space with pointwise operations and the norm $\|\cdot\|_2$ defined above (which arises from an inner product). 
The {\em Hardy algebra} $H^\infty({\mathbb{C}}_+)$ is the set of all  $F\in {\mathcal{O}}({\mathbb{C}}_+)$ such that 
$$
\textstyle
\|F\|_\infty:=\sup\limits_{x>0}  |F(x+iy)|<\infty.
$$
With pointwise operations and the $\|\cdot\|_\infty$ norm, $H^\infty({\mathbb{C}}_+)$ is a Banach algebra. 
With the action of the ring $H^\infty({\mathbb{C}}_+)$ on $H^2({\mathbb{C}}_+)$ given by pointwise multiplication, $H^2({\mathbb{C}}_+)$ is a $H^\infty({\mathbb{C}}_+)$-module. It was shown in \cite[Prop. 8]{Qua} that $H^2({\mathbb{C}}_+)$ is a flat $H^\infty({\mathbb{C}}_+)$-module.

Define $\varphi: {\mathbb{C}}_+\to {\mathbb{D}}$ by 
 \begin{equation}
 \label{30_1_2026_1808}
\textstyle 
\varphi(s)=\frac{s-1}{s+1}\;\text{ for all }\;s\in {\mathbb{C}}_+.
\end{equation}
 Then $\varphi$ is biholomorphic, and
  $
 \textstyle
 \varphi^{-1}(z)=\frac{1+z}{1-z}\;\text{ for all }\;z\in {\mathbb{D}}.
 $ 
It is clear that $f\in H^\infty$ if and only if $f\circ \varphi \in H^\infty({\mathbb{C}}_+)$. 
Equivalently $F\in H^\infty({\mathbb{C}}_+)$ if and only if $F\circ \varphi^{-1} \in H^\infty$. 
By \cite[Theorem, p. 130]{Hof}, $f\in H^2$ if and only if the function 
 $
\textstyle {\mathbb{C}}_+\owns s\mapsto \frac{(f\circ \varphi)(s)}{1+s}
$  
belongs to $ H^2({\mathbb{C}}_+)$. Equivalently, $F\in H^2({\mathbb{C}}_+)$ if and only if 
$$
\textstyle {\mathbb{D}}\owns z\mapsto \frac{(F\circ \varphi^{-1})(z)}{1-z}
$$
 belongs to $H^2$. 
 
\begin{corollary}
Let $\mathfrak{m}$ be a maximal ideal of $H^\infty$ as in Theorem~\ref{sarason}. 
Define $\mathfrak{M}=\{g\circ \varphi: g\in \mathfrak{m}\},$ where $\varphi$ is given by $\eqref{30_1_2026_1808}$. 
 Then $\mathfrak{M}$ is a maximal ideal of $H^\infty({\mathbb{C}}_+),$ and $\mathfrak{M}H^2({\mathbb{C}}_+)=H^2({\mathbb{C}}_+)$. 
\end{corollary}

\noindent {\em Proof}. Let $F\in H^2({\mathbb{C}}_+)$. Then $f\in H^2$, where 
$$
\textstyle 
f(z):=2 \frac{(F\circ \varphi^{-1})(z)}{1-z} \;\text{ for all }\;z\in {\mathbb{D}}. 
$$
Thus there exist $g\in \mathfrak{m}$ and $h\in H^2$ such that $f=gh$. But then $G:=g\circ \varphi\in \mathfrak{M}$, and $H\in H^2({\mathbb{C}}_+)$, where 
$$
 \textstyle H(s)=\frac{(h\circ \varphi)(s)}{1+s} \;\text{ for all }\;s\in {\mathbb{C}}_+.
 $$
 Also, from $f=gh$, we get $f\circ \varphi= (g\circ \varphi)(h\circ \varphi)$, and so for all 
  $s\in {\mathbb{C}}_+$
 $$
 \textstyle 
 \frac{(f\circ \varphi)(s)}{1+s} = (g\circ \varphi)(s) \frac{(h\circ \varphi^{-1})(s)}{1+s}=G(s)H(s).
$$
Thus 
 $
\textstyle 
G(s)H(s)
\!=\! \frac{(f\circ \varphi)(s)}{1+s}\!=\! \frac{1}{1+s} f(\varphi(s))\!=\! \frac{1}{1+s}2\frac{(F\circ \varphi^{-1})(\varphi(s))}{1-\varphi(s)}\!=\! F(s)$. 
$\hfill\Box$ 

\medskip
 
\noindent It can also be shown that if $\mathfrak{N}$ is any finitely generated, nonzero, closed, proper ideal of $H^\infty({\mathbb{C}}_+)$, then $\mathfrak{n}H^2({\mathbb{C}}_+)\subsetneq H^2({\mathbb{C}}_+)$. Indeed, if we set $\mathfrak{n}=\{G\circ \varphi^{-1}: G\in \mathfrak{N}\}$, then $\mathfrak{n}$ is a finitely generated, nonzero, closed, proper ideal of $H^\infty$, and so there exists an $f\in H^2 \setminus (\mathfrak{n}H^2)$. We claim that $F\in H^2({\mathbb{C}}_+)$ defined by 
$$
\textstyle F(s)=2\frac{(f\circ \varphi)(s)}{1+s}\;\text{ for all }\; s\in {\mathbb{C}}_+,
$$
does not belong to $\mathfrak{N} H^2({\mathbb{C}}_+)$. Otherwise, there exist elements $G\in \mathfrak{N} $ and $H\in H^2({\mathbb{C}}_+)$ such that $F=GH$. Define $h\in H^2$ by 
$$
\textstyle
h(z)=\frac{(H\circ \varphi^{-1})(z)}{1-z}\;\text{ for all }\;z\in {\mathbb{D}}.
$$
Thus $F\circ \varphi^{-1}=(G\circ \varphi^{-1})(H\circ \varphi^{-1})$, and so 
 for all $z\in {\mathbb{D}}$, we have 
$$
\textstyle 
\frac{(F\circ \varphi^{-1})(z)}{1-z}=(G\circ \varphi^{-1})(z) \frac{(H\circ \varphi^{-1})(z)}{1-z}=
g(z)h(z).
$$
Thus 
$\textstyle g(z)h(z)\!=\!\frac{(F\circ \varphi^{-1})(z)}{1-z}\!=\!\frac{1}{1-z}F(\varphi^{-1}(z))
\!=\!\frac{1}{1-z}2\frac{(f\circ \varphi)(\varphi^{-1}(z))}{1+\varphi^{-1}(z)}
\!=\!f(z)$. So $f=gh\in \mathfrak{n}H^2$, a contradiction. 
 
\vspace{0.21cm}

\noindent {\bf Acknowledgments.}
 The author thanks Raymond Mortini for simplifying the argument of the strict inclusion $\mathfrak{n}H^2 \!\subsetneq \!H^2$ in \S\ref{Mortini_subsec}, and  Alban Quadrat  for several useful comments, including raising the question of what happens when the maximal ideal $\mathfrak{m}$ under consideration is finitely generated, and for sharing Donald Sarason's email reply.

\end{document}